\documentclass[12pt]{article}
\usepackage{times}
\usepackage[T1]{fontenc}
\usepackage[utf8]{inputenc}
\usepackage[english]{babel}
\usepackage{geometry}
\usepackage{setspace}
\usepackage{amsmath,amssymb,amsthm}
\usepackage{graphicx}
\usepackage{hyperref}
\usepackage{mathtools}
\usepackage{xcolor}
\usepackage{graphicx}
\usepackage{lineno}
\usepackage{rotating}
\usepackage{theoremref}
\usepackage{enumerate}
\usepackage{amstext}
\usepackage[pdf]{graphviz}
\usepackage{mathrsfs}
\newtheorem{theorem}{Theorem}[section]
\newtheorem{lemma}[theorem]{Lemma}
\newtheorem{corollary}[theorem]{Corollary}
\newtheorem{proposition}[theorem]{Proposition}
\theoremstyle{definition}
\newtheorem{definition}[theorem]{Definition}
\newtheorem{example}[theorem]{Example}

\theoremstyle{remark}

\numberwithin{equation}{section}

\title{\textbf{Engel's Theorem for alternative superalgebras}}
\author{Isabel Hern\'andez \\  {\small Secretar\'{\i}a de Ciencia,  Humanidades, Tecnolog\'{\i}a e Innovaci\'on and }\\ {\small Centro de Investigaci\'on en Matem\'aticas  Unidad M\'erida, Mexico.}
\\
\\
Laiz Valim da Rocha \thanks{Email: laizvalim@gmail.com} \\ {\small Instituto de Matem\'atica e Estat\'{\i}stica, Universidade de S\~ao Paulo,  }\\{\small 
S\~ao Paulo, SP, Brazil.}
\\
\\
Rodrigo Lucas Rodrigues\\ {\small Departamento de Matem\'atica, Universidade Federal do Cear\'a, } \\{\small Campus do Pici, Bloco 914, Fortaleza, CE, Brazil.}}
\date{}
\begin{document}
\maketitle
\begin{abstract}
In this paper, a nilpotency criterion is given for finite dimensional alternative superalgebras in the spirit of Engel’s Theorem for Jordan superalgebras over infinite fields provided by Shestakov and Okunev. For alternative superalgebras, no restrictions on the cardinality of the ground field are required. Furthermore, we establish some connections between the concepts of graded-nil and nilpotent alternative superalgebras, and we also exhibit an example of an Engelian commutative power-associative superalgebra of dimension $4$ which is not nilpotent.
\end{abstract}



\section{Introduction}
In the study of nilpotency in nonassociative algebras with associative powers, it is reasonable to ask whether
the nilpotent condition on each element implies
the nilpotency of an algebra. It is well-known that a finite dimensional nil alternative algebra is nilpotent and the same occurs with Jordan algebras over a field of characteristic different from $2$. In finite dimension, the concepts of nil, nilpotent and solvable algebras are equivalent in both classes (see \cite{Schafer}).

In the context of superalgebras, the above result is no longer true. In \cite{superalgebrasandcounter}, I. Shestakov constructed examples of finite dimensional Jordan and alternative superalgebras which are solvable but not nilpotent. Nevertheless, we can observe a relevant aspect in such examples: while the obtained Jordan superalgebra is graded-nil, that is, every homogeneous element is nilpotent, the alternative superalgebra is not. Let us treat some consequences of these facts with the aim of providing nilpotency criteria, starting by considering Jordan superalgebras.

Let $\mathbb{F}$ be a field of characteristic different from $2$. It is immediate to check that a finite dimensional solvable Jordan superalgebra  over $\mathbb{F}$ is graded-nil, since every odd element $x$  satisfies $x^{2}=0$ and the even part is a finite dimensional solvable Jordan algebra, and therefore it is nil.  Consequently, we can conclude that the graded-nil condition is insufficient to guarantee nilpotency for Jordan superalgebras. In fact, this is an expected result, since any odd element $x$ in a Jordan superalgebra over $\mathbb{F}$ is nilpotent, and therefore the graded-nil condition is required only for the even elements. A more strong and appropriate hypothesis
is to impose the nilpotency of the right multiplication operator $R_{a}$, for all homogeneous elements $a\in J$, instead of the nilpotency of $a$, in the sense of Engel's Theorem for Lie algebras \cite{LieJacob}. The latter, fundamental to the development of Lie theory, has been generalized to Lie superalgebras \cite{MR0486011} and several other classes of algebras and superalgebras (for example, see \cite{Engelmalcevsup}, \cite{EngelLeibniz}, \cite{EngelGeneralized} and \cite{Engelforaclass}).

In 2006, K. Okunev and I. Shestakov proved in \cite{EngelJordan} that a finite dimensional Jordan superalgebra $J$ over an infinite field $\mathbb{F}$ of characteristic different from $2$ is nilpotent if it is \emph{Engelian}, that is, the right multiplication operator $R_{a}$ is nilpotent for every homogeneous element $a \in J$. As indicated in \cite[Section 1.1.3]{MR4426392}, the question remains open for Jordan superalgebras over finite fields. 

As mentioned previously, the example of a solvable but not nilpotent finite dimensional alternative superalgebra obtained by Shestakov in \cite{superalgebrasandcounter} is not graded-nil. In fact, unlike what happens with Jordan superalgebras, the odd elements in this case are not necessarily nilpotent. Then we can adapt the initial question for finite dimensional alternative superalgebras and ask if the nil condition on each homogeneous element ensures that the superalgebra is nilpotent. If not, we want to know if it is possible to obtain a nilpotency criterion similar to the Jordan case.

 The purpose of this work is to answer the two questions left above. We exhibit a graded-nil alternative superalgebra over a field of characteristic $3$, which is not nilpotent, in Example \ref{examplenilalternative}. We prove that a finite dimensional alternative superalgebra $A$ is nilpotent if the right multiplication operators by homogeneous elements of $A$ are nilpotent, in Theorem \ref{Engelalternativo}. In particular, a finite dimensional Engelian alternative superalgebra is nilpotent. As shown in Section 2, the notion of Engelian superalgebra is more general than the one presented above for Jordan superalgebras. Finally, we discuss the Engel condition for commutative power-associative superalgebras, presenting an example of a non-Jordan superalgebra of this class of dimension $4$, which is Engelian but not nilpotent. 

\section{Preliminaries}
An algebra is called \emph{alternative} if it satisfies the identities:
$(x,x,y)=0$ and $(x,y,y)=0$, where $(x,y,z)=(xy)z-x(yz)$ denotes the associator of the elements
$x,y,z$. 

A superalgebra $A=A_{0}\oplus A_{1}$ over a field $\mathbb{F}$ is a $\mathbb{Z}_{2}$ graded algebra, that is, $A_{0}$ and $A_{1}$ are vector subspaces of $A$ that satisfy $A_{i}\cdot A_{j}\subseteq A_{i+j \ \text{mod} \ 2}$. A non-zero element $x \in A_{0}\cup A_{1}$ is said to be homogeneous and we denote by $|x|$ the parity index of $x$:
\[ |x|=
\begin{cases}
0, & \text{if } x \in A_{0}\backslash\{0\} \text{ ($x$ is even)};\\
1, & \text{if } x \in A_{1}\backslash\{0\} \text{ ($x$ is odd)}.
\end{cases}
\] 

Let $G={\rm alg}\langle 1, e_{i}, i\in \mathbb{Z}^{+} \ | \ e_{i}e_{j}=-e_{j}e_{i}\rangle$ be the Grassmann algebra $G$ over a field $\mathbb{F}$ . It is easy to see that it has a structure of superalgebra with the grading $G=G_{0}\oplus G_{1}$, where $G_{0}$ is spanned by 1 and the products of an even number of $e_{i}$'s, and $G_{1}$ is spanned by the products of an odd number of $e_{i}$'s. For a superalgebra $A=A_{0}\oplus A_{1}$, we define the \emph{Grassmann envelope} of $A$ as follows: $G(A)=G_{0}\otimes A_{0} \oplus G_{1}\otimes A_{1} $, where the multiplication is 
given by $(x\otimes u)(y\otimes v)=xy\otimes uv,$ 
for every $x \otimes u \in G_{i}\otimes A_{i}, y\otimes v \in G_{j}\otimes A_{j}$, where $i,j \in \{0,1\}$.

\begin{definition}\label{defsup} A superalgebra $A=A_{0}\oplus A_{1}$ is called an \emph{alternative superalgebra} if its Grassman envelope $G(A)$ is an alternative algebra.
\end{definition}

In the same manner, we can define other classes of superalgebras, such as Jordan, Lie and power-associative.
It is important to point out that an associative superalgebra is just an associative algebra with a $\mathbb{Z}_{2}$-graduation and it is also an alternative superalgebra. However, an alternative superalgebra is not always an alternative algebra, as can be seen in \cite{MR1657313}.

Let $A$ be a superalgebra. From now on, we adopt the convention that whenever the parity index appeared in a formula, the corresponding elements are supposed to be homogeneous. By $R_{x}$ and $L_{x}$ we, respectively, denote the right and left multiplication operators by
\[R_{x}(y)=(-1)^{|x||y|}yx \ \ {\rm and} \ \ \ L_{x}(y)=xy.\]

Given a subsuperspace $B$ of $A$, that is, $B=(B\cap A_{0})\oplus (B\cap A_{1})$, we denote by ${B}^{\ast}_{s}$ the (associative) subalgebra of ${\rm End}_{\mathbb{F}} A$ generated by the operators $R_{b}$ and $L_{b}$, where $b$ is homogeneous in $B$. Following the same approach as for the proof of \cite[Theorem 2.4]{Schafer}, we get the following useful result.

\begin{proposition}\label{nilpotenteestrela}
Let $B$ be a subsuperalgebra of a superalgebra $A$. If $B^{\ast}_{s}$ is nilpotent, then $B$ is nilpotent.
\end{proposition}

 By a subsuperalgebra we mean a subsuperspace that is itself a subalgebra. The last proposition tells us that we can prove the nilpotency of a subsuperalgebra $B$ by showing the nilpotency of $B^{\ast}_{s}$, which is associative and, in general, easier to deal with products of its elements.

A homogeneous element $x\in A$ is said to be \emph{Engelian} if the subalgebra of $A^{\ast}_{s}$ generated by $R_{x}$ and $L_{x}$ is nilpotent. A superalgebra $A$ is said to be \emph{Engelian} if every homogeneous element is Engelian. In particular, if $A$ is a supercommutative superalgebra, the Engelian condition is equivalent to $R_{x}$ being nilpotent for every homogeneous element $x\in A$. 
 
\section{A nilpotency criterion for alternative superalgebras} 

From Definition \ref{defsup}, an alternative superalgebra $A=A_{0}\oplus A_{1}$ over a field $\mathbb{F}$ satisfies the following superidentities:
\begin{align}
  (x,y,z)&+(-1)^{|y||z|}(x,z,y)=0\label{lei1}\\ 
  (x,y,z)&+(-1)^{|x||y|}(y,x,z)=0\label{lei2}\\
  (a,a,x)&=0\label{lei3}
\end{align}
where $a \in A_{0}$ and $x,y,z \in (A_{0}\cup A_{1})\backslash\{0\}$.
It is a simple matter to see that the last superidentity follows from the first two in the case of $\rm{char}(\mathbb{F})\neq2$.

By \eqref{lei1} and \eqref{lei2}, we obtain the following superidentities in terms of right and left 
multiplication operators:
\begin{align}
  R_{z}R_{y}&=R_{zy}+(-1)^{|y||z|}\left(R_{yz}-R_{y}R_{z}\right),\label{rr}\\
 L_{z}L_{y}&=L_{zy}+(-1)^{|y||z|}\left(L_{yz}-L_{y}L_{z}\right),\label{ll} \\
 L_{y}R_{z}&= (-1)^{|y||z|}R_{z}L_{y} +L_{yz}-L_{y}L_{z}, \label{lrl}\\
 L_{y}R_{z}&=(-1)^{|y||z|}(R_{z}L_{y} + R_{z}R_{y})-R_{yz}.\label{lrr}
\end{align}

We now give an example of fundamental importance, which will enable us to conclude that a finite dimensional graded-nil alternative superalgebra may not be nilpotent.

\begin{example}\label{examplenilalternative} Let $\mathbb{F}$ be a field of characteristc $3$. Consider the superalgebra $A=A_{0}\oplus A_{1}$ where $A_{0}=\mathbb{F}e_1$ and $A_{1}=\mathbb{F}f_1+\mathbb{F}f_2$, whose non-zero products are $e_{1}f_{2}=f_{2}e_{1}=f_1$ and $f_{1}f_{2}=-f_{2}f_{1}=e_1$.
\end{example}

It is easily seen that $A$ is a nonassociative alternative superalgebra which is graded-nil of index $2$, but is not nilpotent. Moreover, we highlight that considering the superalgebra $A$ over a field of characteristic different from 2, the superalgebra with the same multiplication table is a Jordan superalgebra, as shown in \cite[Example 1]{superalgebrasandcounter}. Consequently, it is convenient to require the Engelian hypothesis in order to obtain a condition of nilpotency for alternative superalgebras. Surprisingly, even if an altenative superalgebra is not necessarily (anti)-commutative, we get a nilpotency criterion just by imposing that the right multiplication operators by homogeneous elements be nilpotent, similarly to Engel's Theorem for Lie algebras.

Throughout this section, let $A=A_{0}\oplus A_{1}$ denote a finite dimensional alternative superalgebra over a field $\mathbb{F}$ such that $R_{a}$ is nilpotent for every homogeneous element $a \in A$. We shall prove, under such conditions, that $A$ is nilpotent. To this end, we follow a similar strategy as in the Jordan case, inspired by \cite{EngelJordan} and \cite[Theorem 3, p. 91]{Nearlyasso}. The next result will be extremely useful to prove our main result.

  \begin{proposition}\label{propnilmaior}
  Let $B$ be a subsuperalgebra of $A$ of the form $B=I+\mathbb{F}v$, where $I$ is a subsuperalgebra of $A$, such that $B^{2}\subseteq I$ and $v \in (A_{0}\cup A_{1})\backslash\{0\}$. The nilpotency of $I_{s}^{\ast}$ implies the nilpotency of $B_{s}^{\ast}$.
 \end{proposition}

To prove the last proposition, we first introduce some notation and auxiliary results. We will denote by $Q$ the right ideal of $B_{s}^{\ast}$ generated by $I_{s}^{\ast}$. It is clear that $Q=I_{s}^{\ast}+I_{s}^{\ast}B_{s}^{\ast}$. It suffices to check the existence of a positive integer $n$ such that $(B_{s}^{\ast})^{n}\subseteq Q$ to conclude the nilpotency of $B_{s}^{\ast}$. In fact, if the previous inclusion occurs for some $n$, we proceed by induction on $m$ that $(B_{s}^{\ast})^{nm}\subseteq (I_{s}^{\ast})^{m}+(I_{s}^{\ast})^{m}B_{s}^{\ast}$ holds for any positive integer $m$. The statement is clear for $m=1$. From the associativity of $B_{s}^{\ast}$, and using the induction hypothesis we have
\begin{align*}
(B_{s}^{\ast})^{n(m+1)}&=(B_{s}^{\ast})^{nm}(B_{s}^{\ast})^{n}\\ 
& \subseteq ((I_{s}^{\ast})^{m}+(I_{s}^{\ast})^{m}B_{s}^{\ast})(B_{s}^{\ast})^{n}\\ 
& \subseteq (I_{s}^{\ast})^{m}(B_{s}^{\ast})^{n}+(I_{s}^{\ast})^{m}B_{s}^{\ast}(B_{s}^{\ast})^{n}\\ 
& \subseteq (I_{s}^{\ast})^{m}(B_{s}^{\ast})^{n} + (I_{s}^{\ast})^{m}(B_{s}^{\ast})^{n+1}\\ 
& \subseteq (I_{s}^{\ast})^{m}(B_{s}^{\ast})^{n} + (I_{s}^{\ast})^{m}(B_{s}^{\ast})^{n} \ \ \text{(since} \ (B_{s}^{\ast})^{n+1}\subseteq (B_{s}^{\ast})^{n} \ \text{)}\\ 
&\subseteq (I_{s}^{\ast})^{m}(I_{s}^{\ast}+I_{s}^{\ast}B_{s}^{\ast}) + (I_{s}^{\ast})^{m}(I_{s}^{\ast}+I^{\ast}B_{s}^{\ast})\\ 
&\subseteq (I_{s}^{\ast})^{m+1} + (I_{s}^{\ast})^{m+1}B_{s}^{\ast}.
\end{align*} 
Hence, the nilpotency of $I_{s}^{\ast}$ implies that of $B_{s}^{\ast}$. 

On the other hand, here and subsequently, we denote by $w=S_{b_{1}}\dots S_{b_{k}}$ a word in $B_{s}^{\ast}$, where $S_{b_{i}}$ is either $R_{b_{i}}$ or $L_{b_{i}}$, and we will call the number $d(w)=k$ the \emph{length} of the word $w$. In addition, since $B=I+\mathbb{F}v$, to prove that $w \in Q$, we can assume, without loss of generality, that each $b_{i}$ is a homogeneous element in $I$ or equal to $v$.

For the purpose of proving the existence of such $n$, the following lemmas will be needed. 
\begin{lemma}\label{aparecea}
Let $w=S_{b_{1}}\dots S_{b_{n}}\in (B_{s}^{\ast})$. If $b_{i}\in I$, for some $1\leq i\leq n$, then $w \in Q$.
\end{lemma}
\begin{proof} Let $a=b_{i}$ given by the hypothesis. If $i=1$, then 
$w \in Q$ follows immediately. Otherwise, we can write
$w=w_{0}S_{a}w_{1}$, where $w_{0}$ is of the form $w_{0}=S_{v}\dots S_{v}$. We will prove the result by induction on the length of $w_{0}$. If $d(w_{0})=1$, then $w_{0}=S_{v}$ 
yields $w=S_{v}S_{a}w_{1}$. We will check that $S_{v}S_{a} \in Q$. To this end, we have the following possibilities:
\vspace{0.2cm}

(I) If $S_{v}=R_{v}$ and $S_{a}=R_{a}$, then
$R_{v}R_{a}=R_{va}\pm(R_{av}-R_{a}R_{v}) \in Q,$ by \eqref{rr};

(II) If $S_{v}=L_{v}$ and $S_{a}=L_{a}$, then
$L_{v}L_{a}=L_{va} \pm(L_{av}-L_{a}L_{v})\in Q,$ by \eqref{ll};

(III) If $S_{v}=L_{v}$ and $S_{a}=R_{a}$, then
$L_{v}R_{a}=\pm (R_{a}L_{v} + R_{a}R_{v})-R_{av} \in Q,$ by \eqref{lrr};

(IV) If $S_{v}=R_{v}$ and $S_{a}=L_{a}$, then
$R_{v}L_{a}=-R_{v}R_{a}+\pm(R_{av}+L_{a}R_{v}) \in Q,$ by \eqref{lrr} and (I).

Now $w \in Q$, which is due to the fact that $Q$ is a right ideal. Assuming $d(w_{0})>1$, we can write $w_{0}=S_{v}w^{\prime}_{0}$, where $w^{\prime}_{0}$ is a subword of $w_{0}$ with $1\leq d(w^{\prime}_{0})< d(w_{0})$. It follows that $w=S_{v}w^{\prime}$, where $w^{\prime}=w^{\prime}_{0}S_{a}w_{1}$. By the induction hypothesis, we have $w^{\prime}\in Q$. Thus, $w^{\prime}$ is a finite sum of words of the form $S_{a_{j}}w_{j}$, with $a_{j} \in I$. Hence, $w$ can be written as a finite sum of words of the form $S_{v}S_{a_{j}}w_{j}$. By the case of length one, we can conclude that each summand belongs to $Q$, so $w \in Q$. In particular,  $Q$ is a two-sided ideal of $B_{s}^{\ast}$. 
\end{proof}

\begin{lemma}\label{aparecesov}
Let $w=S_{v}\dots S_{v} \in (B_{s}^{\ast})^{2n-1}$, where $n$ is an even positive integer such that $R_{v}^{n}=0$. Then $w \in Q$. 
\end{lemma}
 \begin{proof}  
Since $Q$ is a two-sided ideal of $B_{s}^{\ast}$, we can consider the quotient algebra ${B_{s}^{\ast}}/{Q}$. From the superidentities \eqref{lrl}  and \eqref{lrr},  we get 
$- L_{v}L_{v}+L_{v^{2}}=\pm R_{v}R_{v}-R_{v^{2}}.$
In this way, since $R_{v^{2}}+L_{v^{2}}\in I_{s}^{\ast}\subseteq Q$, we obtain $\pm \overline{R_{v}R_{v}}=- \overline{L_{v}L_{v}}.$ Now, $R_{v}^{n}=0$ implies $(\overline{L_{v}L_{v}})^{\frac{n}{2}}= \pm(\overline{R_{v}R_{v}})^{\frac{n}{2}}=\pm\overline{R_{v}^{n}}=\overline{0}$, whence it follows that $L_{v}^{n} \in Q$.

Henceforth, the proof differs according to the parity of $v$.

a) Suppose $v$ is even. As a consequence of \eqref{lei1} and \eqref{lei3}, $R_{v}$ and $L_{v}$ commute. Therefore, we can rewrite $w$ with all left multiplication operators $L_{v}$ in the first positions and then the right multiplication operators $R_{v}$, that is, $w=L_{v}L_{v}\dots R_{v}R_{v}$. Since $d(w)\geq 2n-1$, we deduce that either $R_{v}$ or $L_{v}$ appears at least $n$ times as factors in $w$. In either case, $w \in Q$, since $R_{v}^{n}=0$ and $L_{v}^{n}\in Q$.

b) Assume that $v$ is odd. First, we suppose that $R_{v}$ appears at least $n$ times in $w$. Then \begin{equation}\label{quasicom}
    \overline{L_{v}R_{v}}=-\overline{R_{v}L_{v}}-\overline{R_{v}R_{v}},
\end{equation}
by \eqref{lrr}. We use this relation in ${B_{s}^{\ast}}/{Q}$ to prove that $w \in Q$.

We claim that if $w$ has a subword formed only by the factor $R_{v}$ whose length is $k$, with $k<n$, it is possible, after a finite number of steps, to rewrite 
$\overline{w}$ as a linear combination of words such that each one has a subword formed only by the factor $R_{v}$ whose length is at least $k+1$. Indeed, by hypothesis, $w$ is a word of one of the following forms:
\[w=w_{1}\underbrace{R_{v}\dots R_{v}}_{k \ \text{times}}\underbrace{L_{v}\dots L_{v}}_{m \ \text{times}}R_{v}w_{2} \ \ \text{or} \ \  w=w_{1}R_{v}\underbrace{L_{v}\dots L_{v}}_{m \ \text{times}}\underbrace{R_{v}\dots R_{v}}_{k \ \text{times}}w_{2},\ \ \]
where $w_{1}$ and $w_{2}$ can be empty words. We now proceed 
by induction on $m$. Without lost of generality we can assume
$w=w_{1}\underbrace{R_{v}\dots R_{v}}_{k \ \text{times}}\underbrace{L_{v}\dots L_{v}}_{m \ \text{times}}R_{v}w_{2}$. If $m=1$, then $w=w_{1}\underbrace{R_{v}\dots R_{v}}_{k \ \text{times}}L_{v}R_{v}w_{2}$, and so
 \[
 \begin{aligned}
  \overline{w}&=-\overline{w_{1}\underbrace{R_{v}\dots R_{v}}_{k \ \text{times}}R_{v}L_{v}w_{2}}-\overline{w_{1}\underbrace{R_{v}\dots R_{v}}_{k \ \text{times}}R_{v}R_{v}w_{2}}\\
 &=-\overline{w_{1}\underbrace{R_{v}\dots R_{v}}_{k +1 \ \text{times}}L_{v}w_{2}}-\overline{w_{1}\underbrace{R_{v}\dots R_{v}}_{k+2 \ \text{times}}w_{2}},
 \end{aligned}
 \]
by \eqref{quasicom}, which proves our claim for the case $m=1$.

Suppose that it is satisfied for an arbitrary positive integer $m$. We verify its validity for $m+1$. In this case, note that 
$w=w_{1}\underbrace{R_{v}\dots R_{v}}_{k \ \text{times}}\underbrace{L_{v}\dots L_{v}L_{v}}_{m+1 \ \text{times}}R_{v}w_{2}.$
Application of \eqref{quasicom} gives us
\[
\begin{array}{rcc}
 \overline{w}=&-\overline{w_{1}\underbrace{R_{v}\dots R_{v}}_{k \ \text{times}}\underbrace{L_{v}\dots L_{v}}_{m \ \text{times}}R_{v}L_{v}w_{2}}&-\overline{w_{1}\underbrace{R_{v}\dots R_{v}}_{k \ \text{times}}\underbrace{L_{v}\dots L_{v}}_{m \ \text{times}}R_{v}R_{v}w_{2}} \vspace{0.15cm} \\ 
     & \text{(i)} & \text{(ii)}\\
\end{array}\]

By the induction hypothesis, both (i) and (ii) can be rewritten, after a finite number of steps, as a linear combination of words such that
which one has a subword only formed by the factor $R_{v}$ appearing at least $k+1$ times. Considering all terms obtained, we still have a linear combination and, therefore, with a finite number of steps,  we can rewrite $\overline{w}$ as a linear combination of words such that each one has a subword only formed by the factor $R_{v}$ whose length is at least $k+1$, and the claim is proved.

Applying it to $\overline{w}$, it is possible to write $\overline{w}$ as a linear combination of words such that each one has a subword only formed by the factor $R_{v}$ appearing at least $n$ times, that is, each word is of the form: $\overline{w^{\prime}\underbrace{R_{v}\dots R_{v}}_{n \ \text{times}}w^{\prime \prime}}=\overline{0},$ which implies $\overline{w}=\overline{0}$, that is, $w \in Q$.  Now, if $R_{v}$ does not appear in the writing of $w$ at least $n$ times, so $L_{v}$ appears. By \eqref{lrl}, we infer that
$\overline{L_{v}R_{v}}=-\overline{R_{v}L_{v}}-\overline{L_{v}L_{v}}$

Similarly by changing the roles of $R_{v}$ and $L_{v}$ in our claim, we get that $w \in Q$, and the lemma follows.
\end{proof}

Now, we are in a position to prove Proposition \ref{propnilmaior}.

\begin{proof}\textit{(Proposition \ref{propnilmaior})}
Since $R_{a}$ is nilpotent for every homogeneous element $a\in A$, there exists an even positive integer $n$ such that $R_{v}^{n}=0$. We check that $({B}_{s}^{\ast})^{2n-1}\subseteq Q$ and, consequently, $B_{s}^{\ast}$ will be nilpotent. 

In fact, consider $w=S_{b_{1}}\dots S_{b_{k}}\in (B_{s}^{\ast})^{2n-1}$, with $k\geq 2n-1$, where $b_{i} \in B_{0}\cup B_{1}$, for all
$i=1,\dots,k$. If some $b_{i} \in I$ so, by Lemma \ref{aparecea}, $w \in Q$. Otherwise, $w$ is of the form $w=S_{v_{1}}\dots S_{v_{k}}$, where each $S_{v_{i}}$ is $L_{v}$ or $R_{v}$. In this case, it follows from Lemma \ref{aparecesov} that $w \in Q$. Since each element of $(B_{s}^{\ast})^{2n-1}$ is written as a finite sum of elements of the form of $w$, it also follows that $({B}_{s}^{\ast})^{2n-1}\subseteq Q$, which completes the proof. 
\end{proof}

\begin{theorem}\label{Engelalternativo}
Let $A$ be a finite dimensional alternative superalgebra $A=A_{0}\oplus A_{1}$ over a field $\mathbb{F}$. If $R_{a}$ is nilpotent, for every $a \in A_{0}\cup A_{1}$, then $A$ is nilpotent. In particular, a finite dimensional Engelian alternative superalgebra over a field $\mathbb{F}$ is nilpotent. 
\end{theorem}

\begin{proof}

Note that, if $A=\{0\}$, the result is clear. Therefore, we will assume $A\neq \{0\}$. Let $B$ be a subsuperalgebra strictly contained in $A$ such that $B_{s}^{\ast}$ is nilpotent (a subsuperalgebra with this property always exists, for example, $B=\{0\}$). We prove that it is possible to construct a subsuperalgebra $C=B+\mathbb{F}v$ of $A$, where $v \in (A_{0}\cup A_{1})\backslash\{0\}$ and $v\notin B$, such that $C_{s}^{\ast}$ is also nilpotent. If $C=A$, we have $C_{s}^{\ast}=A_{s}^{\ast}$, which implies the nilpotency of $A$ by Proposition \ref{nilpotenteestrela}. Otherwise, $C$ is a strictly contained subsuperalgebra of $A$ such that $C_{s}^{\ast}$ is nilpotent, and we can repeat this process. Since $A$ is a finite dimensional superalgebra, this process must terminate after finitely many steps. The construction, together with successive applications of Proposition 3.2, yields a strictly increasing sequence of nilpotent subsuperalgebras, leading to the nilpotency of $A$.

\noindent{\bf Claim}. There exists an element $a \in A_{0}\cup A_{1}$ such that $a \notin B$ but $aB\subseteq B$ and $Ba\subseteq B$. 

Suppose that it is false. In this case, for every $a\in A_{0}\cup A_{1}$ such that $a \notin B$, we can find a homogeneous element $b \in B$ such that $ab\notin B$ or $ba \notin B$, i.e., $S_{b}(a)\notin B$, where $S_{b}=R_{b}$ or $S_{b}=L_{b}$. Since $B\neq A$, there exists an element $w_{0}\in A_{0}\cup A_{1}$ such that $w_{0}\notin B$. By the hypothesis, it is possible to find a homogeneous element $b_{1}\in B$ such that $w_{1}=S_{b_{1}}(w_{0})\notin B$. Now, 
since $w_{1}$ is in $A_{0}\cup A_{1}$ but not in $B$, we can find a homogeneous element $b_{2} \in B$ such that $w_{2}=S_{b_{2}}(w_{1})=S_{b_{2}}S_{b_{1}}(w_{0})\notin B$. Let $n$ be the nilpotency index of $B_{s}^{\ast}$. Repeating the procedure $n+1$ times, we obtain a homogeneous element $w_{n+1}\notin B$. On the other hand,    $w_{n+1}=S_{b_{n+1}}\ldots S_{b_{2}}S_{b_{1}}(w_{0})=0$
contradicts $w_{n+1}\notin B$. Hence, there must be $a \in A_{0}\cup A_{1}$ such that $a \notin B$ but $aB\subseteq B$ and $aB\subseteq B$, so the claim is proved.
Now, we have the following possibilities:

(1) Case $a \in A_{0}$. By induction on $k$ and the superidentities \eqref{lei1} and \eqref{lei2}, we can show that $a^{k}B\subseteq B$ and $Ba^{k}\subseteq B$, for every positive integer $k$. Now, since $R_{a}$ is nilpotent and $a \in A_{0}$, so $a$ is nilpotent. Hence, there exists a positive integer $k\geq 2$ such that $a^{k-1}\notin B$ but $a^{k}\in B$. Putting $v=a^{k-1}$, observe that $C=B+\mathbb{F}v$ is a subsuperalgebra of $A$ satisfying $C^{2} \subseteq B$. Thus, by Proposition \ref{propnilmaior}, we can conclude that $C_{s}^{\ast}$ is nilpotent.

(2) Case $a \in A_{1}$. We first assume that $a^{2}\in B$. Here, note that $C=B+\mathbb{F}a$ is a subsuperalgebra satisfying the hypothesis of Proposition \ref{propnilmaior}, hence $C_{s}^{\ast}$ is nilpotent. Now, assume that $a^{2}\notin B$. As in  
(1), we can verify that $a^{2}B\subseteq B$ and $Ba^{2}\subseteq B$. Since $a^{2}\in A_{0}$, follows by (1) that it is possible to find a homogeneous element $v$ such that $C=B+\mathbb{F}v$ is a subsuperalgebra such that $C^{2}\subseteq B$. From the Proposition \ref{propnilmaior}, we obtain that $C_{s}^{\ast}$ is nilpotent.

Therefore, in any of the cases, we can construct a subsuperalgebra $C$ bigger than $B$ such that  $C_{s}^{\ast}$ is nilpotent, and the result is proved. 
\end{proof}
Example \ref{examplenilalternative} exhibits a non-nilpotent   
graded-nil alternative superalgebra of index $2$ 
over a field of characteristic $3$. The next result
confirms that characteristic 3 cannot be removed. 

\begin{corollary}\label{cornil2nilpotent}
Let $A=A_{0}\oplus A_{1}$ be a graded-nil finite dimensional alternative superalgebra over a field $\mathbb{F}$ of characteristic different from $3$. If every $a_{1}\in A_{1}$ is nilpotent of index $2$,
then $A$ is nilpotent.
\end{corollary}
\begin{proof}
First, observe that since every odd element is nilpotent of index $2$, if $a,b\in A_{1}$, then $ab=-ba$. Now, suppose that $A$ is not nilpotent. By Theorem \ref{Engelalternativo}, there exists a homogeneous element $f$ such that $R_{f}$ is not nilpotent. Note that $f$ must be an odd element. In fact, since $R_{e}^{k}=R_{e^{k}}$, for every even element $e\in A_{0}$, and $A$ is graded-nil, it follows that $e$ is nilpotent, and $R_{e}^{k}$ is as well. Therefore, $f\in A_{1}$. 

Now,  $R_{f}$ being non-nilpotent implies $R_{f}R_{f}\neq 0$. Accordingly, there exists an even element $e$ such that $R_{f}R_{f}(e)\neq 0$, otherwise $R_{f}^{3}=0$, and, in consequence, $(e,f,f)=(ef)f-ef^{2}=-R_{f}R_{f}(e)\neq 0.$
On the other hand, \eqref{lei2} gives $(ef)f=-(fe)f+f(ef)$. Consequently, $2f(ef)=(fe)f$. Also, from \eqref{lei1}, $(fe)f-f(ef)=f(fe)$, and hence $2(fe)f=f(ef)$. If $\rm{char}(\mathbb{F})=2$, we obtain $(e,f,f)=(ef)f=-f(ef)=0$, which is absurd. If $\rm{char}(\mathbb{F})\neq 2$, we get $4f(ef)=2(fe)f=f(ef)$.  Hence, $3f(ef)=0$. In this case, we also obtain $(e,f,f)=(ef)f=-f(ef)=0$, which again is a contradiction.
\end{proof}

Finally, we also observe that the superalgebra of Example \ref{examplenilalternative} is nonassociative. In fact, for associative superalgebras, a homogeneous element $x$ is nilpotent if and only if $R_{x}$ is nilpotent, since these superalgebras are associative algebras. Thus, it follows from Theorem \ref{Engelalternativo}:

\begin{corollary}
A finite dimensional graded-nil associative superalgebra is nilpotent.    
\end{corollary}
\section{Engelian commutative power-associative superalgebras}

An algebra is said to be power-associative if every subalgebra generated by a single element is associative.
It is well known that alternative and Jordan algebras belong to the class of power-associative algebras, see \cite{Schafer}. In particular, Jordan algebras are contained in the commutative subclass and satisfy $(x,y,x^2)=0$. Thus, a natural question was posed by Albert \cite{Albert}: as in the case of Jordan algebras in characteristic different from $2$, is every finite dimensional nil commutative power-associative algebra also nilpotent? The answer is negative. Suttles provided in \cite{suttles} a counterexample of dimension $5$, while Gerstenhaber and Myung showed in \cite{Gerstenhaber} that the result holds up to dimension $4$ over a field of characteristic different from $2$.

Over a field of characteristic prime to $30$, the variety of commutative power-associative superalgebras can be described by polynomial identities, namely, the commutative law and $x^{2}x^{2}=(x^{2}x)x$. In the context of superalgebras, the corresponding superidentities can be found in \cite{lowcommu}. Due to the supercommutativity, as in the Jordan case, it is suitable to replace nilpotency with the Engelian condition on homogeneous elements and consider the following question: Is every finite dimensional Engelian commutative power-associative superalgebra nilpotent? Once again, the answer is negative. Indeed, the example provided by Suttles has a trivial superalgebra structure that is Engelian and not nilpotent. But we can ask whether this is the best possible dimension, as in the case of commutative-power associative algebras.

In order to answer this question, we can consider the list of all commutative power-associative superalgebras up to dimension $4$ over an algebraically closed field $\mathbb{F}$, with $\mathrm{char}(\mathbb{F})$ prime to $30$, provided by Hern\'andez et al. in 2021 \cite{lowcommu}. The superalgebra $(2,2)_{57}^{\ast}$ described in that paper is a non-Jordan example of Engelian superalgebra which is not nilpotent. 

In fact, with the homogeneous basis $\{e_1,e_2,f_1,f_2\}$, with $|e_i|=0$ and $|f_j|=1$ for any $1\leq i,j\leq 2$, the only nonzero products of basis elements in $(2,2)_{57}^{\ast}$ are:
\begin{equation*}
e_1^{2}=e_2, \ e_2f_2=f_2e_2=f_1, \ f_1f_2=-f_2f_1=e_1. 
\end{equation*}

It is simple to check that $R_{x}$ is nilpotent with the nilindex $2$, if $x$ is even and with the nilindex $3$, if $x$ is odd in $(2,2)_{57}^{\ast}$. Moreover, it satisfies $((f_1 f_2)e_1)f_2 = f_1$, which implies the non-nilpotency of this superalgebra.

In dimension $\leq 3$ all Engelian superalgebras founded in \cite{lowcommu} are Jordan, so it is not possible to find counterexamples in lower dimensions in this case.

\section*{Acknowledgments}
The authors wish to express their gratitude to Prof. Ivan Shestakov for sending us a detailed email with suggestions for an appropriate approach to the problem presented in this paper. The first author was
supported by Grant CONAHCYT A1-S-45886. This study was financed in part by the Coordenação de Aperfeiçoamento de Pessoal de Nível Superior – Brasil (CAPES) – Finance Code 001.

\end{document}